\documentclass[11pt]{article}
\usepackage{amsmath,amssymb}
\usepackage{amsfonts}
\usepackage[dvipsnames]{xcolor}
\usepackage[utf8]{inputenc}
\usepackage[T1]{fontenc}
\usepackage{float} \usepackage[colorlinks=true]{hyperref}
\hypersetup{urlcolor=blue, citecolor=blue, linkcolor=blue}
\newtheorem{theo}{Theorem}[section]
\newtheorem{lem}[theo]{Lemma}
\newtheorem{cor}[theo]{Corollary}
\newtheorem{prop}[theo]{Proposition}

\newtheorem{defi}[theo]{Definition}

\newcommand{\mysection}[1]{\section{#1} \setcounter{equation}{0}}
\newcommand{\proof}{{\sc Proof.} \quad}
\newcommand{\proofc}{{\sc Proof} \ }
\newcommand{\be}{\begin{equation} \label}
\newcommand{\ee}{\end{equation}}
\newcommand{\bea}{\begin{eqnarray}\label}
\newcommand{\eea}{\end{eqnarray}}
\newcommand{\bas}{\begin{eqnarray*}}
\newcommand{\eas}{\end{eqnarray*}}
\newcommand{\bit}{\begin{itemize}}
\newcommand{\eit}{\end{itemize}}
\newcommand{\qed}{\hfill$\Box$ \vskip.2cm}
\newcommand{\nn}{\nonumber}
\newcommand{\R}{\mathbb{R}}
\newcommand{\N}{\mathbb{N}}
\newcommand{\pO}{\partial\Omega}

\newcommand{\eps}{\varepsilon}

\newcommand{\supp}{{\rm supp} \, }

\newcommand{\wto}{\rightharpoonup}

\newcommand{\hra}{\hookrightarrow}
\newcommand{\io}{\int_\Omega}

\newcommand{\bom}{\overline{\Omega}}

\newcommand{\abs}{\\[5pt]}

\newcommand{\uie}{u_{i\eps}}
\newcommand{\uene}{u_{1\eps},...,u_{N\eps}}
\newcommand{\iin}{i\in \{1,...,N\}}
\newcommand{\jin}{j\in \{1,...,N\}}
\newcommand{\fiu}{f_i(u_1,...,u_N)}
\newcommand{\fiue}{f_i(u_{1\eps},...,u_{N\eps})}
\newcommand{\fjue}{f_j(u_{1\eps},...,u_{N\eps})}
\newcommand{\weps}{w_{\eps}}
\newcommand{\zeps}{z_{\eps}}
\newcommand{\oze}{\overline{z_{\eps}}}

\newcommand{\calU}{\mathcal{U}}

\newcommand{\reactrates}[2]{\overset{#1}{\underset{#2}{\rightleftharpoons}}}

\usepackage{newunicodechar}
\newunicodechar{α}{\alpha}
\newunicodechar{β}{\beta}
\newunicodechar{γ}{\gamma}
\newunicodechar{λ}{\lambda}
\newunicodechar{ε}{\varepsilon}

\begin{document}
\enlargethispage{10mm}
\title{
Global existence in reaction-diffusion systems with mass control under relaxed assumptions
merely referring to cross-absorptive effects}
\author{
Johannes Lankeit\footnote{lankeit@ifam.uni-hannover.de}\\
{\small Leibniz Universität Hannover, Institut f\"ur Angewandte Mathematik,}\\
{\small Welfengarten~1, 30167 Hannover, Germany}
\and
Michael Winkler\footnote{michael.winkler@math.uni-paderborn.de}\\
{\small Institut f\"ur Mathematik, Universit\"at Paderborn,  }\\
{\small Warburger Str.~100, 33098 Paderborn, Germany}
}
\date{}
\maketitle
\begin{abstract}
\noindent
We introduce a generalized concept of solutions for reaction-diffusion systems and prove their global existence. 
The only restriction on the reaction function beyond regularity, quasipositivity and mass control is special in that it merely  controls the growth of cross-absorptive terms.\\
The result covers nonlinear diffusion and does not rely on an entropy estimate.\abs
\end{abstract}
\newpage
\mysection{Introduction}\label{intro}
Reaction-diffusion equations arise in various applications in chemistry and biology (cf. \cite[Ch. 2]{pierre_survey} and form an important class of model problems in the study of systems of parabolic equations (see \cite[Ch. 33]{quittner_souplet}).
Already at the stage of basic theories of solvability, a major challenge for the analysis of such systems 
consists in the presence of commonly superlinear source terms.
While the possibility of blow up then is apparrent as long as suitably destabilizing reaction mechanisms are admitted (cf. e.g. \cite{laamri_91}), even the requirement of dissipation of mass -- which is sufficient to yield global existence and boundedness in the corresponding ODE systems -- cannot preclude its occurrence, as impressively demonstrated by the counterexamples in \cite{pierre_schmitt_blowup}. 
Global classical solutions have, accordingly, been searched for and found under certain restrictive conditions. \abs
In the context of boundary value problems for systems of the general form
\be{00}
	\partial_t u_i = d_i \Delta u_i^{m_i} 
	+ \fiu,
	\qquad i\in \{1,...,N\},
\ee
in the linear diffusion case when $m_1=...=m_N=1$ such results on global smooth solvability
cover settings where boundedness of the first among two components is a priori known (from a sign of $f_1$, cf.~\cite{martin_pierre97}),
where the diffusion coefficients are close to each other (\cite{cupps_morgan_tang,canizo_desvillettes_fellner}), or where sources
exhibit subquadratic growth (\cite{caputo_vasseur}), 
and recently ideas of \cite{kanel} have successfully been extended to show global solvability 
for quadratic or slightly superquadratic reaction functions (\cite{souplet18,caputo_goudon_vasseur,fellner_morgan_tang}). \abs
Another line of investigations pursues solutions in a weaker sense. 
Weak solutions can be constructed if $L^1$-bounds for the reaction terms $f_i(u)$ are known (\cite{pierre03}), or if the reaction functions are at most quadratic (\cite{pierre_rolland}). For nonlinear diffusion of porous medium type, the existence of weak solutions to \eqref{0} with Dirichlet boundary data is shown in \cite{laamri_pierre} under the 
assumptions that
\be{f2}
	\left\{ \begin{array}{l}
	f_i\in C^1([0,\infty)^N)
	\quad \mbox{is such that} \quad \\[1mm]
	f_i(s_1,...s_{i-1},0,s_{i+1},...,s_N) \ge 0
	\qquad \mbox{for all $(s_1,...,s_N)\in [0,\infty)^N$ and } \iin,
	\end{array} \right.
\ee
that with some $K \ge 0$ and $a\in(0,\infty)^N$
we have
\be{f1}
	\sum_{i=1}^N a_i f_i(s_1,...,s_N) \le K \cdot \bigg( \sum_{i=1}^N a_i  s_i+1\bigg)
	\qquad \mbox{for all $(s_1,...,s_N)\in [0,\infty)^N$ and } \iin,
\ee
and that
either a priori $L^1$-bounds for the reaction terms are known, or $f_i(s_1,..,s_N) \ge - C \left(\sum_{j=}^n s_j^{β_j} +1 \right)$ with $β_i<m_i+1$ for all $i\in\{1,\ldots,N\}$ 
An analogous result has been achieved for the corresponding Neumann problem in \cite{laamri_perthame} with a different proof and less restrictive conditions on the initial data.\abs
%
The concept of renormalized solutions, that is, the idea that not $u$ itself, but a transformed quantity $\rho(u)$ solves (a weak form) of the equation, makes it possible to bypass even further restrictions on the form of the system. This concept has been successfully introduced for the Boltzmann equation by DiPerna and Lions \cite{diperna_lions} and was employed for reaction diffusion equations with quadratic reaction functions and linear diffusion in \cite{desvillettes_fellner_pierre_vovelle}. The apparently most far-reaching application of this idea to reaction-diffusion systems (with linear diffusion) can be found in \cite{fischer_ARMA}, where essentially no growth restriction on the $f_i$ is needed, but where the reaction function is supposed to obey a certain entropy condition. The term in the definition of solutions for whose treatment this entropy condition is essential arises from the choice of renormalization functions $\xi\colon [0,\infty)^N\to \R$ with compactly supported $D\xi$, which in particular depend on all solution components simultaneously.\abs
{\bf Main results.} \quad
In the present manuscript, we intend to introduce an approach by which it becomes possible to avoid any requirement of the latter
type, and it turns out that this can in fact be achieved by resorting to separate renormalization functions for each component $u_i$. Thereby, our main result, as stated in Theorem~\ref{theo9} below, 
partially answers the open problem \cite[p.585]{fischer_ARMA} to find a similar notion of solution without requiring an entropy condition.\abs
Specifically, we shall be concerned with the Neumann problem
\be{0}
	\left\{
	\begin{array}{ll}
	\partial_t u_i = d_i \Delta u_i^{m_i} 
	+ \fiu,
	\qquad & x\in\Omega, \ t>0, \ \iin,\\[1mm]
	\partial_\nu u_i^{m_i}=0,
	\qquad & x\in\pO, \ t>0, \ \iin,\\[1mm]
	u_i(x,0)=u_{0i}(x),
	\qquad & x\in\Omega, \ \iin,
	\end{array} \right.
\ee
under the assumptions that (\ref{f2}) and (\ref{f1}) hold, and that
\be{f3}
	\left\{	\begin{array}{l}
	f_i(s_1,..,s_N) \ge - \phi_i(s_i) \cdot \bigg( 
	\sum\limits_{\begin{array}{c} 
	\\[-7mm]
	\scriptscriptstyle
	j\in\{1,...,N\} \\[-2.5mm]
	\scriptscriptstyle
	j\ne i
	\end{array}}
	s_j^{\beta_j}
	+ 1 \bigg)
	\quad \mbox{for all $(s_1,...,s_N)\in [0,\infty)^N$ and } \iin, \\
	\mbox{with some nonnegative $\phi_i\in C^1([0,\infty))$ such that $\phi_i'>0$ on $(0,\infty)$ for $\iin$, and some} 
	\\[1mm]
	\beta_i>0 \mbox{ such that $\beta_i<m_i+1$ for all $\iin$.}
	\end{array} \right.
\ee
Here we recall that the quasipositivity condition in \eqref{f2} is important in order to avoid negative concentrations, and that
\eqref{f1} is a slightly generalized mass dissipation condition, and includes some stoichiometric coefficients $a$.
In addition to this, \eqref{f3} signifies a growth condition for the negative parts of the reaction functions, where in the special case of linear diffusion, subquadratic growth is admissible. It is important to note, however, that this restriction only applies to the cross-absorptive effects: For $(f_i)_-$, the possible growth with respect to the $i$-th argument remains unrestricted.\abs
As for the initial data in (\ref{0}), throughout this paper we shall suppose that
\be{init}
	\mbox{$u_{0i}$, $\iin$, is a nonnegative function from $L^r(\Omega)$ with some }
	\left\{ \begin{array}{ll}
	r\ge 1 \quad & \mbox{if } n=1, \\[1mm]
	r>1 \quad & \mbox{if } n=2, \\[1mm]
	r\ge \frac{2n}{n+2} \quad & \mbox{if } n\ge 3.
	\end{array} \right.
\ee
Postponing the precise description of the solution concept to be pursued here to Section~\ref{sec:solutionconcept},
let us introduce our main result obtained in this framework, and give a few examples of its application. 
\begin{theo}\label{theo9}
  let $n\ge 1$ and $N\ge 1$ and $\Omega \subset\R^n$ a bounded domain with smooth boundary, and suppose that $d_1,...,d_N$ and
  $m_1,...,m_N$ are positive, and that $f_1,...,f_N$ belong to $C^1([0,\infty))$ and satisfy (\ref{f1}), (\ref{f2})
  and (\ref{f3}) with some positive constant $K$.
  Then given any $u_{01},...,u_{0N}$ fulfilling (\ref{init}), one can find nonnegative functions 
  $u_i\in L^{m_i+1}_{loc}(\bom\times [0,\infty))$ such that $(u_1,...,u_N)$ is a generalized solution of (\ref{0})
  in the sense of Definition \ref{defi6}.
\end{theo}
\vspace*{6mm}
{\bf Application \#1.} \quad		
%
%
%
%
%
%
%
A first application of Theorem \ref{theo9} addresses the system
\be{11.3}
	\left\{
	\begin{array}{ll}
	\partial_t u_i = d_i \Delta u_i^{m_i} 
	+ (p_i-q_i) \cdot \bigg( k_2 \prod_{j=1}^N u_j^{q_j} - k_1 \prod_{j=1}^N u_j^{p_j}\bigg),
	\qquad & x\in\Omega, \ t>0, \ \iin,\\[1mm]
	\partial_\nu u_i^{m_i}=0,
	\qquad & x\in\pO, \ t>0, \ \iin,\\[1mm]
	u_i(x,0)=u_{0i}(x),
	\qquad & x\in\Omega, \ \iin,
	\end{array} \right.
\ee
which describes a general reversible reaction of the form 
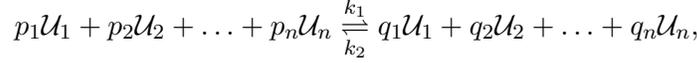
\[
 p_1\calU_1+p_2\calU_2+\ldots+p_n\calU_n \reactrates{k_1}{k_2} q_1\calU_1+q_2\calU_2+\ldots+q_n\calU_n,
\]
and for which we obtain the following.
\begin{prop}\label{prop11}
  Let $N\ge 2$, and suppose that $k_1>0$ and $k_2>0$, and that for $i\in\{1,...,N\}$, $d_i>0$, $m_i>0$, $p_i\ge 1$ and 
  $q_i\ge 1$ are such that for some $a\in (0,\infty)^N$
  \be{stoch}
	\sum_{i=1}^N a_i p_i=\sum_{i=1}^N a_i q_i,
  \ee
  and that
  \be{11.1}
	\sum_{\begin{array}{c} 
	\\[-6mm]
	\scriptstyle
	j\in\{1,...,N\} \\[-1mm]
	\scriptstyle
	j\ne i
	\end{array}}
	\frac{p_j}{m_j+1} <1
	\qquad \mbox{for all $\iin$ such that } p_i>q_i
  \ee
  as well as
  \be{11.2}
	\sum_{\begin{array}{c} 
	\\[-6mm]
	\scriptstyle
	j\in\{1,...,N\} \\[-1mm]
	\scriptstyle
	j\ne i
	\end{array}}
	\frac{q_j}{m_j+1} <1
	\qquad \mbox{for all $\iin$ such that } p_i<q_i.
  \ee
  Then for any choice of $u_{01},...,u_{0n}$ complying with (\ref{init}), the problem (\ref{11.3})
  admits a generalized solution in the sense of Definition \ref{defi6}.
\end{prop}
\proofc of Proposition \ref{prop11}. \quad
  Writing $f_i(s_1,...,s_N):=(p_i-q_i) \cdot \Big(k_2 \prod_{j=1}^N s_j^{q_j} - k_1 \prod_{j=1}^N s_j^{p_j}\Big)$
  for $\iin$ and $(s_1,...,s_N) \in [0,\infty)^N$, we see that \eqref{f2} is fulfilled and \eqref{f1} follows since 
  $\sum_{j=1}^N a_i  f_i \equiv 0$ due to (\ref{stoch}).
  Moreover, if e.g.~$\iin$ is such that $p_i>q_i$, then (\ref{11.1}) enables us to pick numbers $\theta_j>1$,
  $\jin \setminus \{i\}$, such that $p_j\theta_j <m_j+1$ for all $\jin\setminus \{i\}$ and 
  $\sum_{j\ne i} \frac{1}{\theta_j} <1$.
  An application of Young's inequality thus shows that for any such $i$,
  \bas
	f_i(s_1,...,s_N)
	&\ge& -(p_i-q_i) k_1 s_i^{p_i} \cdot \prod_{j\ne i} s_j^{p_j} \\
	&\ge& -(p_i-q_i) k_1 s_i^{p_i} \cdot \bigg( \sum_{j\ne i} s_j^{p_j \theta_j} +1 \bigg)
	\qquad \mbox{for all } (s_1,...,s_N) \in [0,\infty)^N,
  \eas
  and complementing this by a similar reasoning for all $\iin$ for which $p_i<q_i$, 
  we readily obtain that holds and Theorem \ref{theo9}
  becomes applicable so as to yield the claim.
\qed
Proposition \ref{prop11} corresponds to \cite[Remark 2.10]{laamri_pierre}, where the existence of weak solutions is proved. The main difference is that there the summation in \eqref{11.1} and \eqref{11.2} extends over all $j\in\{1,\ldots N\}$.\abs
%
For linear diffusion, weak solutions of \eqref{11.3} have been found in \cite{pierre_suzuki_umakoshi} if the reaction functions grow at most quadratically or if the diffusion coefficients are sufficiently close to each other. The same article also deals with their exponential convergence.

\vspace*{4mm}
{\bf Application \#2.} \quad		
A second application of our general theory is concerned with the variant of (\ref{11.3}) given by
\be{ex2}
	\left\{
	\begin{array}{ll}
	\partial_t u_1 = d_1 \Delta u_1^{m_1} 
	+ u_1^{\beta_1} g_2(u_2) - g_1(u_1) u_2^{\beta_2},
	\qquad & x\in\Omega, \ t>0, \\[1mm]
	\partial_t u_2 = d_2 \Delta u_2^{m_2} 
	- u_1^{\beta_1} g_2(u_2) + \lambda g_1(u_1) u_2^{\beta_2},
	\qquad & x\in\Omega, \ t>0, \\[1mm]
	\partial_\nu u_1^{m_1}=\partial_\nu u_2^{m_2}=0,
	\qquad & x\in\pO, \ t>0, \\[1mm]
	u_1(x,0)=u_{01}(x),
	\quad u_2(x,0)=u_{02}(x),
	\qquad & x\in\Omega,
	\end{array} \right.
\ee
and underlines the mildness of the assumptions in Theorem \ref{theo9} by admitting widely arbitrary growth of the main ingredients
$g_1$ and $g_2$ appearing herein:
\begin{prop}\label{prop13}
  Let $d_1>0, d_2>0, m_1>0$, $m_2>0$ and $\lambda\in [0,1]$, let $\beta_1\in [1,m_1+1)$ and $\beta_2\in [1,m_2+1)$, and let
  $g_1\in C^1([0,\infty))$ and $g_2\in C^1([0,\infty))$ be such that $g_1(0)=g_2(0)=0$ and that
  $g_1$ and $g_2$ are positive on $(0,\infty)$.
  Then for any pair $(u_{01}, u_{02})$ satisfying (\ref{init}), there exists a generalized solution of (\ref{ex2})
  in the spirit of Definition \ref{defi6}.
\end{prop}
\proofc of Proposition \ref{prop13}. \quad
  Taking any nonnegative $\phi_i\in C^1([0,\infty))$ such that $\phi_i'>0$ and $\phi_i \ge g_i$ on $(0,\infty)$ for 
  $i\in\{1,2\}$, one can readily verify that for
  \bas
	f_1(s_1,s_2):=s_1^{\beta_1} g_2(s_2) - g_1(s_1) s_2^{\beta_2}
	\quad \mbox{and} \quad
	f_2(s_1,s_2):=-s_1^{\beta_1} g_2(s_2) +\lambda g_1(s_1) s_2^{\beta_2},
	\qquad (s_1,s_2)\in [0,\infty)^2,
  \eas
  we have
  \bas
	f_1(s_1,s_2)+f_2(s_1,s_2)= -(1-\lambda) g_1(s_1) s_2^{\beta_2} \le 0
	\qquad \mbox{for all } (s_1,s_2)\in [0,\infty)^2
  \eas
  as well as
  \bas
	f_1(s_1,s_2) \ge -g_1(s_1) s_2^{\beta_2}
	\ge -\phi_1(s_1) s_2^{\beta_2}
	\qquad \mbox{for all } (s_1,s_2)\in [0,\infty)^2
  \eas
  and, similarly,
  \bas
	f_s(s_1,s_2) \ge -\phi_2(s_2) s_1^{\beta_1}
	\qquad \mbox{for all } (s_1,s_2)\in [0,\infty)^2.
  \eas
  The assumptions $\beta_i<m_i+1$, $i\in \{1,2\}$, therefore warrant applicability of Theorem \ref{theo9} with the 
  intended result.
\qed
Let us remark that
since in Proposition~\ref{prop13} not only $f_1+f_2\le 0$, but also $\lambda f_1(s_1,s_2)+f_2(s_1,s_2)=(\lambda-1) s_1^{\beta_1}g_2(s_2)\le 0$ for $(s_1,s_2)\in[0,\infty)^2$, \cite[Cor. 2.11]{laamri_pierre} could be applied to the variant of \eqref{ex2} 
involving homogeneous Dirichlet boundary conditions (cf. \cite[Remark 2.12]{laamri_pierre}) so as to yield weak solutions for any $L^1$-initial data; said corollary, however, requires that $m_1,m_2<2$.

\vspace*{4mm}
{\bf Application \#3.} \quad		
We shall next briefly address
\be{ex3}
	\left\{
	\begin{array}{ll}
	\partial_t u_1 = d_1 \Delta u_1^{m_1} 
	+ k_2 u_1^{q_1} u_2^{q_2} - k_1 u_1^{p_1} u_2^{p_2},
	\qquad & x\in\Omega, \ t>0, \\[1mm]
	\partial_t u_2 = d_2 \Delta u_2^{m_2} 
	- k_2 u_1^{q_1} u_2^{q_2} + k_1 u_1^{p_1} u_2^{p_2},
	\qquad & x\in\Omega, \ t>0, \\[1mm]
	\partial_\nu u_1^{m_1}=\partial_\nu u_2^{m_2}=0,
	\qquad & x\in\pO, \ t>0, \\[1mm]
	u_1(x,0)=u_{01}(x),
	\quad u_2(x,0)=u_{02}(x),
	\qquad & x\in\Omega,
	\end{array} \right.
\ee
for which without imposing any smallness condition on $q_2$ nor $p_1$ we obtain the following.
\begin{cor}\label{cor12}
  Let $k_1,k_2,d_1,d_2,m_1$ and $m_2$ be positive, and let $p_1\ge 1, p_2\ge 1, q_1\ge 1$ and $q_2\ge 1$ be such that
  \bas
	q_1<m_1+1
	\qquad \mbox{and} \qquad
	p_2<m_2+1.
  \eas
  Then for all $(u_{01}, u_{02})$ fulfilling (\ref{init}), the problem (\ref{ex3})
  possesses a generalized solution in the sense of Definition \ref{defi6}.
\end{cor}
\proofc of Corollary \ref{cor12}. \quad
 With $g_1(s)=k_1s^{p_1}$, $g_2(s)=k_2s^{q_2}$, $\beta_1=q_1$, $\beta_2=p_2$, this immediately results from Proposition \ref{prop13}. 
\qed

\vspace*{4mm}
{\bf Application \#4.} \quad		
As final example, let us consider the generalized Lotka-Volterra system
\be{LV}
	\left\{
	\begin{array}{ll}
	\partial_t u_i = d_i \Delta u_i^{m_i} 
	+ \gamma_i u_i + \sum_{j=1}^N a_{ij} u_j^{\beta_{ij}}u_i^{\beta_{ji}},
	\qquad & x\in\Omega, \ t>0, \ \iin,\\[1mm]
	\partial_\nu u_i^{m_i}=0,
	\qquad & x\in\pO, \ t>0, \ \iin,\\[1mm]
	u_i(x,0)=u_{0i}(x),
	\qquad & x\in\Omega, \ \iin,
	\end{array} \right.,
\ee
which does not obey the typical entropy condition (that is required for the renormalized solutions in \cite{fischer_ARMA} and for  classical solvability e.g. in \cite{souplet18}). In \cite{fellner_morgan_tang}, global classical solutions are shown to exist for the classical Lotka-Volterra system ($\beta_{ij}=1$ for all $i,j$) with linear diffusion. If $m_i=1$ for all $i$ and $\beta_{ij}=\beta$ for all $i,j$, then the result of \cite{fellner_morgan_tang} covers $2\beta < 2+\varepsilon$ (for sufficiently small $\varepsilon\in(0,1)$, see \cite[Step 1, (7)]{fellner_morgan_tang}).\abs
Within the generalized solvability framework considered here, the following consequence of Theorem \ref{theo9}
shows that here actually the entire range $\beta<2$ can be exhausted.
\begin{prop}\label{prop:LV}
  Let $N\ge 2$, $d_i>0$, $m_i>0$ and $\gamma_i\in \R$ for $i\in \{1,\ldots,N\}$, and suppose that 
  for $i,j\in\{1,...,N\}$ the numbers $a_{ij}\in\R$ and $\beta_{ij}>0$ are such that
  $a_{ij}+a_{ji}\le 0$, and that
  \begin{equation}\label{LVcondaij}
  	\text{if } a_{ij}<0, \quad \text{ then } \beta_{ij}<m_i+1.
  \end{equation}
  Then for all initial data  $u_{01},...,u_{0n}$ as in \eqref{init}, \eqref{LV} 
  has a generalized solution in the sense of Definition \ref{defi6}.
\end{prop}
\newcommand{\sumiN}{\sum_{i=1}^N}
\newcommand{\sumjN}{\sum_{j=1}^N}
\proofc of Proposition \ref{prop:LV}.\quad 
  With $f_i(s)=\gamma_i s_i + \sum_{j=1}^N a_{ij} s_j^{\beta_ij}s_i^{\beta_{ji}}$, $i\in\{1,\ldots,N\}$, \eqref{f2} is clearly   
  satisfied. As 
  \begin{align*}
 	\sum_{i=1}^N f_i(s) &= \sum_{i=1}^N \gamma_i s_i + \sumiN \sumjN a_{ij} s_j^{\beta_{ij}}s_i^{\beta_{ji}}
 	&=\sum_{i=1}^N \gamma_i s_i + \sum_{i\neq j} (a_{ij}s_j^{\beta_{ij}}s_i^{\beta_{ji}}+a_{ji}s_i^{\beta_{ji}}s_j^{\beta_{ij}})
 	&\le \max_i \gamma_i \sumiN s_i,
  \end{align*}
  also \eqref{f1} holds. Finally, 
  \begin{align*}
 	f_i(s) &\ge \gamma_i s_i + \sum\limits_{\begin{array}{c} 
	\\[-7mm]
	\scriptscriptstyle
	j\in\{1,...,N\} \\[-2.5mm]
	\scriptscriptstyle
	a_{ij}<0
	\end{array}} a_{ij} s_j^{\beta_{ij}}s_i^{\beta_{ji}} \ge - \bigg(|\gamma_i| + \sum\limits_{\begin{array}{c} 
	\\[-7mm]
	\scriptscriptstyle
	j\in\{1,...,N\} \\[-2.5mm]
	\scriptscriptstyle
	a_{ij}<0
	\end{array}} |a_{ij}| s_j^{\beta_{ij}} \bigg)\bigg(s_i+\sumjN s_i^{\beta_{ji}}\bigg)\\
	&\ge - \phi_i(s) \bigg(1+\sum\limits_{\begin{array}{c} 
	\\[-7mm]
	\scriptscriptstyle
	j\in\{1,...,N\} \\[-2.5mm]
	\scriptscriptstyle
	a_{ij}<0
	\end{array}} s_j^{\beta_{ij}}\bigg)
  \end{align*}
  for all $s=(s_1,...,s_N)\in[0,\infty)^N$ 
  if we set $\phi(s)=\max\{|\gamma_i|, |a_{ij}|\mid j\in\{1,\ldots N\}\}\cdot \sumjN s_i^{\beta_{ji}}$ for any such $s$, so
  that, according to \eqref{LVcondaij}, \eqref{f3} is fulfilled and Theorem~\ref{theo9} is applicable.
\qed 
\mysection{Solution concept}\label{sec:solutionconcept}
The first step toward the design of our solution concept is concerned with an appropriate supersolution feature required
in each of the equations making up (\ref{0}):
\begin{defi}\label{defi5}
  Suppose that for $\iin$, $u_{0i}:\Omega\to \R$ and $u_i:\Omega\times (0,\infty) \to\R$ are measurable and nonnegative.
  Then $(u_1,...,u_N)$ will be called a {\em renormalized supersolution} of (\ref{0}) if for every $\rho\in C^\infty([0,\infty))$
  fulfilling $\rho'\in C_0^\infty([0,\infty))$, $\rho'\le 0$ and $\rho''\ge 0$, with
  \be{5.1}
	P_i^{(1)}(s):=\int_0^s \sigma^\frac{m_i-1}{2} \sqrt{\rho''(\sigma)} d\sigma
	\quad \mbox{and} \quad
	P_i^{(2)}(s):=\int_0^s \sigma^{m_i-1} \rho'(\sigma) d\sigma,
	\qquad s\ge 0, \ \iin,
  \ee
  we have
  \be{5.23}
	\rho'(u_i) \fiu \in L^1_{loc}(\bom\times [0,\infty))
	\quad \mbox{and} \quad
	\nabla P_i^{(1)}(u_i) \in L^2_{loc}(\bom\times [0,\infty);\R^n)
	\qquad \mbox{for all } \iin,
  \ee
  and if moreover
  \bea{5.4}
	- \int_0^\infty \io \rho(u_i)\varphi_t 
	- \io \rho(u_{0i}) \varphi(\cdot,0)
	&\le& - d_i m_i \int_0^\infty \io \varphi |\nabla P_i^{(1)}(u_i)|^2 
	+ d_1 m_i \int_0^\infty \io P_i^{(2)}(u_i) \Delta \varphi \nn\\
	& & + \int_0^\infty \io \rho'(u_i) \fiu \varphi
  \eea
  for all $\iin$ and each nonnegative $\varphi\in C_0^\infty(\bom\times [0,\infty))$ fulfilling
  $\partial_\nu\varphi=0$ on $\pO\times (0,\infty)$.
\end{defi}
{\bf Remark.} \quad
i) \ In the above situation, both integrals on the left of (\ref{5.4}) as well as the second integral on the 
right-hand side therein exist due to the readily verified fact that
$\rho$ and $P_i^{(2)}$, $\iin$, are bounded on $[0,\infty)$.\abs
ii) \ The supersolution property in \cite[Prop. 3.6]{laamri_pierre} is obtained upon the choice of $\rho(x)=x$ (inadmissible in Definition~\ref{defi5}), integration by parts in the integral involving $P_i^{(2)}$ (and addition of a corresponding integrability requirement) and, finally, exchange of $\leq$ by $\geq$ (resulting from the change of sign of $\rho'$).

\vspace*{4mm}
As discussed in several previous related approaches toward generalized solvability on the basis of
supersolution features of the above flavor (\cite{zhigun_SIMA}, \cite{win_JEE2018}),
supplementing Definition \ref{defi5} by 
a mere requirement on mass control is already sufficient to create a notion of solvability which within classes
of suitably smooth functions indeed reduces to classical ones (see, e.g, \cite{lankeit_win_NoDEA} and \cite{win_SIMA2015}
for detailed reasonings in this regard):
\begin{defi}\label{defi6}
  By a {\em generalized solution} of (\ref{0}) we mean a vector $(u_1,...,u_N)$ of nonnegative measurable functions
  on $\Omega\times (0,\infty)$ such that $(u_1,...,u_N)$ is a renormalized supersolution of (\ref{0}) in the sense of
  Definition \ref{defi5}, that with some $a\in(0,\infty)^N$
  \be{6.1}
	u_1,...,u_N 
	\quad \mbox{and} \quad 
	\sum_{i=1}^N a_i \fiu
	\quad \mbox{belong to } L^1_{loc}(\bom\times [0,\infty)),
  \ee
  and that
  \be{6.2}
	\io \bigg(\sum_{i=1}^N a_i  u_i(\cdot,t)\bigg)
	\le \io \bigg( \sum_{i=1}^N a_i u_{0i} \bigg)
	+ \int_0^t \io \bigg( \sum_{i=1}^N a_i  \fiu \bigg)
	\qquad \mbox{for a.e.~} t>0.
  \ee
\end{defi}
\mysection{Approximate systems}
In order to construct such solutions through an essentially standard type of approximation,
for $\eps\in (0,1)$ we consider
\be{0eps}
	\left\{
	\begin{array}{ll}
	\displaystyle
	\partial_t \uie = d_i \Delta (\uie+\eps)^{m_i} 
	+ \frac{\fiue}{1+\eps\sum_{j=1}^N \fiue},
	\qquad & x\in\Omega, \ t>0, \ \iin,\\[1mm]
	\partial_\nu \uie=0,
	\qquad & x\in\pO, \ t>0, \ \iin,\\[1mm]
	\uie(x,0)=u_{0i\eps}(x),
	\qquad & x\in\Omega, \ \iin,
	\end{array} \right.
\ee
where
\be{ie}
	\left\{ \begin{array}{l}
	(u_{0i\eps})_{\eps\in (0,1)} \subset C^1(\bom)
	\quad \mbox{ is such that $u_{0i\eps} \ge 0$ for all $\iin$, that}\\[1mm]
	u_{0i\eps} \to u_{0i}
	\quad \mbox{in $L^1(\Omega)$ and a.e.~in $\Omega$ as $\eps\searrow 0$ for all $\iin$, and that} \\[1mm]
	\sup_{\eps\in (0,1)} \|u_{0i\eps}\|_{L^r(\Omega)} <\infty
	\end{array} \right.
\ee
with $r\ge 1$ taken from (\ref{init}).\abs
Due to boundedness of the reaction term therein and nondegeneracy of the diffusion, by \cite[Theorems 14.4 and 14.6]{amann} and \cite[Chapter V]{LSU}, 
for each fixed $\eps\in(0,1)$ the problem \eqref{0eps} indeed admits a global classical solution.\abs
\textbf{General assumption.} Throughout the sequel, we shall suppose that the assumptions of Theorem~\ref{theo9} and \eqref{ie} are satisfied, and given $\eps\in(0,1)$ we let $(\uene)$ denote the global classical solution of \eqref{0eps}.\abs
The following basic observation concerning $L^1$-boundedness of these solutions is a fairly immediate consequence of (\ref{f1}).
\begin{lem}\label{lem1}
  For all $T>0$ there exists $C(T)>0$ such that for all $\iin$ and any $\eps\in (0,1)$ we have
  \be{1.1}
	\|\uie(\cdot,t)\|_{L^1(\Omega)} \le C(T)
	\qquad \mbox{for all } t\in (0,T).
  \ee
\end{lem}
\proof
  By integrating in (\ref{0eps}), we see that since $\partial_\nu \uie=0$ on $\pO\times (0,\infty)$,
  due to (\ref{f1}) we have
  \bas
	\frac{d}{dt} \sum_{i=1}^N a_i \io \uie
	&=& \io \frac{1}{1+\eps\sum_{j=1}^N |\fjue|}
	\cdot \sum_{i=1}^N a_i \fiue \\
	&\le& \io \frac{1}{1+\eps\sum_{j=1}^N |\fjue|}
	\cdot K \bigg\{ \sum_{i=1}^N a_i \uie+1\bigg\} \\
	&\le& K \sum_{j=1}^N a_i \io \uie + K|\Omega|
	\qquad \mbox{for all $t>0$ and } \eps\in (0,1),
  \eas
  because $\uie\ge 0$ for all $\iin$. By an ODE comparison, this shows that
  \bas
	\sum_{i=1}^N a_i \io \uie
	\le \bigg\{ \sum_{i=1}^N a_i \io u_{0i\eps} + |\Omega| \bigg\} \cdot e^{Kt}
	\qquad \mbox{for all $t>0$ and } \eps\in (0,1)
  \eas
  and hence, again by nonnegativity of $\uie$ for $\iin$, establishes (\ref{1.1}) due to (\ref{ie}).
\qed
The following estimate rests on a duality-based reasoning inspired by a corresponding argument from \cite{laamri_pierre}. An important difference is given by the change of boundary conditions: Where \cite{laamri_pierre} dealt with Dirichlet boundary data and thus could use the solution of Poisson's equation as test function, the non-invertibility of the Laplacian with homogeneous Neumann boundary data leads us to employ the solution of a Helmholtz equation instead. An alternative approach of working with the Neumann Laplacian after substracting the mean value has been followed in \cite{laamri_perthame}, but it seems unclear how far strategies of this type can be applied so as to successfully cover the present setting.\abs
\begin{lem}\label{lem2}
  For all $T>0$ there exists $C(T)>0$ such that for all $\iin$,
  \be{2.1}
	\int_0^T \io \uie^{m_i+1} \le C(T)
	\qquad \mbox{for all } \eps\in (0,1).
  \ee
\end{lem}
\proof
  Following \cite[Proof of Theorem 2.7]{laamri_pierre}, we first observe that according to (\ref{0eps}) and (\ref{f1}), writing
  \bas
	\weps(\cdot,t):= e^{-Kt} \cdot \sum_{i=1}^N a_i \uie(\cdot,t)
	\qquad \mbox{and} \qquad
	\zeps(\cdot,t):=\int_0^t e^{-Ks} \cdot \bigg\{ \sum_{i=1}^N a_i d_i(\uie+\eps)^{m_i} \bigg\}
  \eas
  for $\eps\in (0,1)$ and $t\ge 0$, we have
  \bas
	\partial_t \weps
	&=& e^{-Kt} \Delta \bigg\{ \sum_{i=1}^N a_i d_i (\uie+\eps)^{m_i} \bigg\}
	+ e^{-Kt} \sum_{i=1}^N \frac{\fiue}{1+\eps\sum_{j=1}^N |\fjue|} 
	- K e^{-Kt} \sum_{i=1}^N a_i \uie \\
	&\le& e^{-Kt} \Delta \bigg\{ \sum_{i=1}^N a_i d_i (\uie+\eps)^{m_i} \bigg\}
	+ K e^{-Kt}
	\quad \mbox{in $\Omega\times (0,\infty)$ \quad for all } \eps\in (0,1)
  \eas
  and hence
  \be{2.2}
	\weps(\cdot,t)
	\le \weps(\cdot,0) + \Delta \zeps(\cdot,t) +1
	\quad \mbox{in $\Omega$ \quad for all $t>0$ and } \eps\in (0,1),
  \ee
  because $\int_0^t K e^{-Ks} ds = 1-e^{-Ks}\le 1$ for all $t>0$.
  Upon multiplication by $\partial_t \zeps \ge 0$ and integration over $\Omega\times (0,T)$ for $T>0$, as in \cite{laamri_pierre}
  we obtain that since $\zeps(\cdot,0)\equiv 0$,
  \bea{2.3}
	\int_0^T \io \weps \partial_t \zeps
	&\le& \io \weps(\cdot,0) \zeps(\cdot,T) 
	- \frac{1}{2} \io |\nabla\zeps(\cdot,T)|^2
	+ \io \zeps(\cdot,T) \nn\\
	&\le& \io \Big(\weps(\cdot,0)+1\Big) \zeps(\cdot,T) 
	\qquad \mbox{for all } \eps\in (0,1),
  \eea
  where the right-hand side is now estimated in a way slightly deviating from that in \cite{laamri_pierre} due to the different boundary
  conditions considered here. 
  In fact, by nonnegativity of $\weps$ the inequality in (\ref{2.2}) implies that for each fixed $T>0$ and arbitrary 
  $\eps\in (0,1)$,
  \bas
	\left\{ \begin{array}{l}
	-\Delta \zeps(\cdot,T) + \zeps(\cdot,T)
	\le \weps(\cdot,0)+1+\zeps(\cdot,T)
	\quad \mbox{in $\Omega$ \qquad and } \\[1mm]
	\partial_\nu \zeps(\cdot,T)=0
	\quad \mbox{on } \pO,
 	\end{array} \right.
  \eas
  so that since the Helmholtz operator $-\Delta+1$ admits a comparison principle under homogeneous Neumann boundary conditions,
  we obtain that 
  \be{2.4}
	\zeps(\cdot,T) \le \oze
	\qquad \mbox{in } \Omega,
  \ee
  where $\oze$ denotes the solution of
  \bas
	\left\{ \begin{array}{l}
	-\Delta \oze + \oze = \weps(\cdot,0)+1+\zeps(\cdot T)
	\quad \mbox{in } \Omega, \\[1mm]
	\partial_\nu \oze=0
	\quad \mbox{on } \pO.
	\end{array} \right.
  \eas
  Now without loss of generality assuming that the number $r\ge 1$ in (\ref{init}) satisfies $r<2$, we take $r'\in (2,\infty]$
  such that $\frac{1}{r}+\frac{1}{r'}=1$, and employ a Sobolev embedding theorem and elliptic regularity theory
  (\cite{GT}) to find $c_1>0$ and $c_2>0$ such that for all $\eps\in (0,1)$,
  \begin{align}\label{2.5}
	\|\oze\|_{L^{r'}(\Omega)}
	& \le c_1\|\oze\|_{W^{2,r}(\Omega)}
	\le c_2\|\weps(\cdot,0)+1+\zeps(\cdot,T)\|_{L^{r}(\Omega)} \nn\\
	& \le c_2 \|\weps(\cdot,0)+1\|_{L^r(\Omega)} + c_2\|\zeps(\cdot,T)\|_{L^r(\Omega)},
  \end{align}
  because the restrictions on $r$ in (\ref{init}) warrant that $2-\frac{n}{r}\ge -\frac{n}{r'}$ if $n\ge 3$, and that
  $2-\frac{n}{r}>0$ if $n\le 2$.
  Here we note that according to Lemma \ref{lem1} we know that if we let $\theta\in (0,1)$ small enough such that
  $\theta m_i \le 1$ for all $\iin$, then we can find $c_3(T)>0$ fulfilling
  \bas
	\io \zeps^\theta(\cdot,T)
	&\le& T^\theta \sup_{t\in(0,T)} \io \bigg\{ \sum_{i=1}^Na_i d_i (\uie(\cdot,t)+\eps)^{m_i} \bigg\}^\theta \\
	&\le& \Big\{ NT \max_{\iin}a_i  d_i \Big\}^\theta \sup_{t\in(0,T)}\io (\uie(\cdot,t)+\eps)^{\theta m_i} \\
	&\le& \Big\{ NT \max_{\iin} a_i d_i \Big\}^\theta \sup_{t\in(0,T)}\io (\uie(\cdot,t)+1) \\
	&\le& c_3(T)
	\qquad \mbox{for all } \eps\in (0,1),
  \eas
  whence using that then $\theta<1<r<2<r'$ we may invoke Young's inequality to see that with some $c_4>0$ we have
  \bas
	c_2\|\zeps(\cdot,T)\|_{L^r(\Omega)}
	&\le& \frac{1}{2} \|\zeps(\cdot,T)\|_{L^{r'}(\Omega)}
	+ c_4\|\zeps(\cdot,T)\|_{L^\theta(\Omega)} \\
	&\le& \frac{1}{2} \|\zeps(\cdot,T)\|_{L^{r'}(\Omega)}
	+ c_3^\frac{1}{\theta}(T) c_4
	\qquad \mbox{for all } \eps\in (0,1).
  \eas
  In view of (\ref{2.4}) and the nonnegativity of $\zeps$, (\ref{2.5}) thus implies that
  \bas
	\|\zeps(\cdot,T)\|_{L^{r'}(\Omega)}
	\le \|\oze\|_{L^{r'}(\Omega)}
	\le c_2 \|\weps(\cdot,0)+1\|_{L^r(\Omega)}
	+ \frac{1}{2}\|\zeps(\cdot,T)\|_{L^{r'}(\Omega)}
	+ c_3^\frac{1}{\theta}(T) c_4
	\qquad \mbox{for all } \eps\in (0,1),
  \eas
  so that in (\ref{2.3}) we can use the H\"older inequality to estimate
  \bas
	\io \Big(\weps(\cdot,0)+1\Big)\zeps(\cdot,T)
	&\le& \|\weps(\cdot,0)+1\|_{L^r(\Omega)} 
	\|\zeps(\cdot,T)\|_{L^{r'}(\Omega)} \\
	&\le& \|\weps(\cdot,0)+1\|_{L^r(\Omega)} 
	\cdot \Big\{
	2c_2 \|\weps(\cdot,0)+1\|_{L^r(\Omega)} + 2c_3^\frac{1}{\theta}(T) c_4 \Big\}
  \eas
  for all $\eps\in (0,1)$.
  As $\sup_{\eps\in (0,1)} \|\weps(\cdot,0)\|_{L^r(\Omega)}$ is finite according to the hypothesis (\ref{ie}),
  (\ref{2.3}) thereby entails the existence of $c_5(T)>0$ such that
  \bas	
	\int_0^T \io e^{-Kt} \cdot \bigg\{ \sum_{i=1}^N d_i(\uie+\eps)^{m_i} \bigg\}
	\cdot \bigg\{ \sum_{i=1}^N \uie \bigg\}
	\le c_5(T)
	\qquad \mbox{for all } \eps\in (0,1),
  \eas
  from which (\ref{2.1}) readily follows.
\qed
We next rely on \eqref{f3} in deriving the following estimates for gradients and reaction terms. Testing \eqref{0eps} by $-\frac{1}{\phi_i(\uie)}$, namely, enables us to successfully combine \eqref{f3} with \eqref{2.1}. In order to obtain a bound for, e.g., $|\fiue|$ instead of $|-\frac{1}{\phi_i(\uie)}\fiue|$, we here restrict our attention to sets of the form $\{\uie\le M\}$, where  $|-\frac{1}{\phi_i(\uie)}|$ can be estimated from below by a positive constant.
\begin{lem}\label{lem3}
  Let $M>0$ and $T>0$. Then one can find $C(M,T)>0$ such that
  \be{3.1}
	\int_0^T \io \chi_{\{\uie\le M\}} (\uie+\eps)^{m_i-1} |\nabla\uie|^2 \le C(M,T)
	\qquad \mbox{for all $i\in\N$ and } \eps\in (0,1)
  \ee
  and
  \be{3.2}
	\int_0^T \io \chi_{\{\uie\le M\}} 
	\frac{|\fiue|}{1+\eps\sum_{j=1}^N |\fjue|}
	\le C(M,T)
	\qquad \mbox{for all $i\in\N$ and } \eps\in (0,1).
  \ee
\end{lem}
\proof
  For fixed $\iin$, we take $\phi_i\in C^1([0,\infty))$ as in (\ref{f3}) and define
  \bas
	\Phi_i(s):=-\int_1^{s+1} \frac{d\sigma}{\phi_i(\sigma)},
	\qquad s\ge 0.
  \eas
  Then since $\phi_i$ is nondecreasing, we have
  \be{3.21}
	0 \ge \Phi_i(s) \ge - c_{1i} s
	\qquad \mbox{for all } s\ge 0
  \ee
  with $c_{1i}:=\frac{1}{\phi_i(1)}>0$, and moreover $\Phi_i'(s)=-\frac{1}{\phi_i(s+1)}, \ s\ge 0$, satisfies
  \be{3.22}
	0 \le - \Phi_i'(s) \le \frac{1}{\phi_i(s)}
	\qquad \mbox{for all } s\ge 0.
  \ee
  Now using (\ref{0eps}), for $\eps\in (0,1)$ and $t>0$ we compute
  \bas
	\frac{d}{dt} \io \Phi_i(\uie)
	&=& \io \Phi_i'(\uie) \cdot \bigg\{ d_i m_i \nabla \cdot \Big((\uie+\eps)^{m_i-1} \nabla\uie\Big)
	+ \frac{\fiue}{1+\eps\sum_{j=1}^N |\fjue|} \bigg\} \\
	&=& - d_i m_i \io (\uie+\eps)^{m_i-1} \Phi_i''(\uie) |\nabla\uie|^2
	+ \io \Phi_i'(\uie) \cdot \frac{\fiue}{1+\eps\sum_{j=1}^N |\fjue|},
  \eas
  so that splitting $f_i=|f_i|-2(f_i)_-$, upon further integration we find that
  \bea{3.3}
	& & \hspace*{-20mm}
	d_i m_i \int_0^T \io (\uie+\eps)^{m_i-1} \Phi_i''(\uie) |\nabla\uie|^2
	+ \int_0^T \io |\Phi_i'(\uie)| \cdot \frac{|\fiue|}{1+\eps\sum_{j=1}^N |\fjue|} \nn\\
	&=& \io \Phi_i(u_{0i\eps}) - \io \Phi_i(\uie(\cdot,T)) \nn\\
	& & + 2\int_0^T \io |\Phi_i'(\uie)| \cdot \frac{(\fiue)_-}{1+\eps\sum_{j=1}^N |\fjue|}
	\qquad \mbox{for all } T>0,
  \eea
  because $\Phi_i'$ is nonpositive.
  Now on the right-hand side of (\ref{3.3}) we use (\ref{f3}) to see that thanks to (\ref{3.22}),
  \bas
	2\int_0^T \io |\Phi_i'(\uie)| \cdot \frac{(\fiue)_-}{1+\eps\sum_{j=1}^N |\fjue|}
	&\le& 2\int_0^T \io |\Phi_i'(\uie)| \cdot \phi_i(\uie) \cdot \bigg\{ \sum_{j\ne i} u_{j\eps}^{\beta_j}+1\bigg\} \\
	&\le& 2\sum_{j\ne i} \int_0^T \io u_{j\eps}^{\beta_j}
	+ 2|\Omega|T
	\qquad \mbox{for all } T>0,
  \eas
  whence recalling Lemma \ref{lem2} we obtain $c_{2i}(T)>0$ such that for all $\eps\in (0,1)$,
  \be{3.4}
	2\int_0^T \io |\Phi_i'(\uie)| \cdot \frac{(\fiue)_-}{1+\eps\sum_{j=1}^N |\fjue|}
	\le c_{2i}(T).
  \ee
  As
  \bas
	\io \Phi_i(u_{0i\eps}) - \io \Phi_i(\uie(\cdot,T))
	\le c_{1i} \io \uie(\cdot,T)
	\qquad \mbox{for all } T>0
  \eas
  by (\ref{3.21}), in view of Lemma \ref{lem1} we thus conclude from (\ref{3.3}) that there exists $c_{3i}(T)>0$ with the
  property that
  \be{3.5}
	d_i m_i \int_0^T \io (\uie+\eps)^{m_i-1} \Phi_i''(\uie) |\nabla\uie|^2
	+ \int_0^T \io |\Phi_i'(\uie)| \cdot \frac{|\fiue|}{1+\eps\sum_{j=1}^N |\fjue|}
	\le c_{3i}(T)
  \ee
  for all $\eps\in (0,1)$.
  Noting that for fixed $M>0$ and any $\iin$ the numbers
  \bas
	c_{4i}:=\min_{s\in [0,M]} \Phi_i''(s)
	= \min_{s\in [0,M]} \frac{\phi_i'(s+1)}{\phi_i^2(s+1)}
	\qquad \mbox{and} \qquad
	c_{5i}:=\min_{s\in [0,M]} |\Phi_i'(s)| = \frac{1}{\phi_i(M+1)}
  \eas
  are both positive according to our hypotheses on $\phi_i$ and $\phi_i'$, from (\ref{3.5}) we readily infer (\ref{3.1})
  and (\ref{3.2}) if we let $C(M,T):=\max\limits\Big\{ \max_{\iin} \frac{c_{3i}(T)}{d_i m_i c_{4i}} \, , \, 
  \max_{\iin} \frac{c_{3i}(T)}{c_{5i}} \Big\}$, for instance.
\qed
While the bound in Lemma~\ref{lem2} is sufficient for concluding relative compactness of $\{\uie\mid \eps\in(0,1)\}$ in some weak topology, we are additionally interested in possible pointwise convergence of $u_{i\eps_j}$ along some sequence $(\eps_j)_{j\in\N}\searrow0$. 
We thus strive to derive a suitable strong compactness property in $L^2(\Omega\times(0,T))$, at least of a power of $\uie$
which has been cut off at large values so as to ensure accessibility to the estimates of Lemma~\ref{lem3}.
\begin{lem}\label{lem4}
  Given $\zeta\in C_0^\infty([0,\infty))$, for $\iin$ let $\rho_i(s):=s^{\kappa_i}\zeta(s)$, $s\ge 0$, where
  $\kappa_i:=\max\{\frac{m_i+1}{2} \, , \, 2\}$.
  Then
  \be{4.1}
	\Big(\rho_i(\uie)\Big)_{\eps\in (0,1)}
	\mbox{ is relatively compact in $L^2(\Omega\times (0,T))$ for all } T>0.
  \ee
\end{lem}
\proof
  Let us first make make sure that for each $T>0$,
  \be{4.2}
	\Big(\rho_i(\uie)\Big)_{\eps\in (0,1)}
	\mbox{ is bounded in } L^2((0,T);W^{1,2}(\Omega)).
  \ee
  To see this, we note that due to the compactness of $\supp \rho_i$ it is clear that 
  $(\rho_i(\uie))_{\eps\in (0,1)}$ is bounded in
  $L^\infty(\Omega\times (0,\infty))$. Thereore, (\ref{4.2}) results upon the observation that 
  if we fix $M>0$ such that $\zeta\equiv 0$ on $(M,\infty)$, then by 
  Young's inequality and (\ref{3.1}) we see that for all $T>0$ there exists $c_1(T)>0$ such that
  \bas
	\int_0^T \io |\nabla\rho_i(\uie)|^2
	&=& \int_0^T \io (\rho_i'(\uie))^2 |\nabla\uie|^2 \\
	&=& \int_0^T \io \Big( \kappa_i \uie^{\kappa_i-1} \zeta(\uie) + \uie^{\kappa_i} \zeta'(\uie)\Big)^2 |\nabla\uie|^2 \\
	&\le& 2\kappa_i^2 \|\zeta\|_{L^\infty((0,\infty))}^2 
		\int_0^T \io \chi_{\{\uie\le M\}} \uie^{2\kappa_i-2} |\nabla\uie|^2 \\
	& & + 2\|\zeta'\|_{L^\infty((0,\infty))}^2 
		\int_0^T \io \chi_{\{\uie\le M\}} \uie^{2\kappa_i} |\nabla\uie|^2 \\
	&\le& 2\kappa_i^2 \|\zeta\|_{L^\infty((0,\infty))}^2 M^{2\kappa_i-m_i-1} 
		\int_0^T \io \chi_{\{\uie\le M\}} \uie^{m_i-1} |\nabla\uie|^2 \\
	& & + 2\|\zeta'\|_{L^\infty((0,\infty))}^2 M^{2\kappa_i-m_i+1}
		\int_0^T \io \chi_{\{\uie\le M\}} \uie^{m_i-1} |\nabla\uie|^2 \\[2mm]
	&\le& c_1(T)
	\qquad \mbox{for all } \eps\in (0,1),
  \eas
  where we have used that $2\kappa_i-m_i+1 \ge 2\kappa_i-m_i-1\ge 0$ by hypothesis.\abs
  We next fix any integer $k\ge 1$ such that $k>n$, and claim that then for all $T>0$,
  \be{4.3}
	\Big(\partial_t \rho_i(\uie)\Big)_{\eps\in (0,1)}
	\mbox{ is bounded in } L^1 \Big( (0,T);(W^{k,2}(\Omega))^\star \Big).
  \ee
  To verify this in quite a straightforward manner, we pick $\psi\in C^\infty(\bom)$ and use (\ref{0eps}) to see that for 
  each $t>0$ and any $\eps\in (0,1)$,
  \bea{4.4}
	& & \hspace*{-20mm}
	\bigg| \io \partial_t \rho_i(\uie(\cdot,t)) \cdot\psi \bigg| \nn\\
	&=& \Bigg| \io \rho_i'(\uie) \psi \cdot 
	\bigg\{ d_i m_i \nabla \cdot \Big( (\uie+\eps)^{m_i-1} \nabla\uie \Big)
	+ \frac{\fiue}{1+\eps\sum_{j=1}^N |\fjue|} \bigg\} \Bigg| \nn\\
	&=& \Bigg| -d_i m_i \io (\uie+\eps)^{m_i-1} \rho_i''(\uie) |\nabla\uie|^2 \psi
	- d_i m_i \io (\uie+\eps)^{m_i-1} \rho_i'(\uie) \nabla\uie\cdot\nabla\psi \nn\\
	& & \hspace*{10mm}
	+ \io \rho_i'(\uie) \cdot \frac{\fiue}{1+\eps\sum_{j=1}^N |\fjue|} \cdot \psi \Bigg| \nn\\
	&\le& d_i m_i \|\rho_i''\|_{L^\infty((0,\infty))} \cdot
	\bigg\{ \io \chi_{\{\uie\le M\}} (\uie+\eps)^{m_i-1} |\nabla\uie|^2 \bigg\} \cdot \|\psi\|_{L^\infty(\Omega)} \nn\\
	& & + d_i m_i \cdot \bigg\{
	\io \chi_{\{\uie\le M\}} (\uie+\eps)^{m_i-1} |\nabla\uie|^2
	+ |\Omega|\sup_{s\ge 0} (s+\eps)^{m_i-1} |\rho_i'(s)|^2 \bigg\} \|\nabla \psi\|_{L^\infty(\Omega)}\nn\\
	& & + \|\rho_i'\|_{L^\infty((0,\infty))} \cdot 
	\bigg\{ \io \chi_{\{\uie\le M\}} \frac{|\fiue|}{1+\eps\sum_{j=1}^N |\fjue|} \bigg\} \cdot \|\psi\|_{L^\infty(\Omega)},
  \eea
  where finiteness of both $\|\rho'\|_{L^\infty((0,\infty))}$ and $\|\rho_i''\|_{L^\infty((0,\infty))}$
  is asserted by our restriction that $\kappa_i\ge 2$.
  Here we observe that in the case $m_i\ge 1$ we have
  \bas
	(s+\eps)^{m_i-1} |\rho_i'(s)|^2 \le (M+1)^{m_i-1} \|\rho_i'\|_{L^\infty((0,\infty))}^2,
  \eas
  whereas if $m_i<1$ then
  \bas
	(s+\eps)^{m_i-1} |\rho_i'(s)|^2 
	&\le& s^{m_i-1} |\rho_i'(s)|^2 \\
	&=& s^{m_i-1} \cdot \Big|\kappa_i s^{\kappa_i-1} \zeta(s) + s^{\kappa_i} \zeta'(s)\Big|^2 \\
	&\le& 2\kappa_i^2 s^{2\kappa_i+m_i-3} \zeta^2(s)
	+ 2s^{2\kappa_i+m_i-1} |\zeta'(s)|^2 \\
	&\le& 2\kappa_i^2 M^{2\kappa_i+m_i-3} \|\zeta\|_{L^\infty((0,\infty))}^2
	+ 2M^{2\kappa_i+m_i-1} \|\zeta'\|_{L^\infty((0,\infty))}^2,
  \eas
  because again by definition of $\kappa_i$, we have $2\kappa_i+m_i-1\ge 2\kappa_i+m_i-3 \ge 0$.
  As furthermore $W^{k,2}(\Omega) \hra W^{1,\infty}(\Omega)$ due to our restriction that $k>n$,
  from (\ref{4.4}) we thus infer the existence of $c_2>0$ such that
  for all $t>0$ and any $\eps\in (0,1)$,
  \bas
	\Big\| \partial_t \rho_i(\uie(\cdot,t))\Big\|_{(W^{k,2}(\Omega))^\star}
	&\le& c_2 \io \chi_{\{\uie\le M\}} (\uie+\eps)^{m_i-1} |\nabla\uie|^2  + c_2 \\
	& & + c_2 \io \chi_{\{\uie\le M\}}  \frac{|\fiue|}{1+\eps\sum_{j=1}^N |\fjue|},
  \eas
  which in light of (\ref{3.1}) and (\ref{3.2}) establishes \eqref{4.3} upon a time integration.\abs
  We finally only need to combine (\ref{4.2}) with (\ref{4.3}) to conclude that (\ref{4.1}) is a consequence of
  an Aubin-Lions type lemma (\cite[Cor. 4]{simon}). 
\qed
Based on this compactness statement, we may conclude the existence of a limit object.
\begin{lem}\label{lem7}
There exist $(\eps_j)_{j\in\N} \subset (0,1)$ and nonnegative functions $u_1,...,u_N$ defined on $\Omega\times (0,\infty)$
  such that $\eps_j\searrow 0$ as $j\to\infty$, and such that for all $\iin$,
  \be{7.1}
	\uie \to u_i
	\qquad \mbox{a.e.~in } \Omega\times (0,\infty)
  \ee
  and
  \be{7.2}
	\uie \to u_i
	\qquad \mbox{in } L^p_{loc}(\bom\times [0,\infty))
	\quad \mbox{for all } p\in [1,m_i+1)
  \ee
  as well as
  \be{7.3}
	\uie \wto u_i
	\qquad \mbox{in } L^{m_i-1}_{loc}(\bom\times [0,\infty))
  \ee
  as $\eps=\eps_j\searrow 0$.
\end{lem}
\proof
  In Lemma \ref{lem4} choosing $\zeta=\zeta_l$, with $\zeta_l \in C_0^\infty([0,\infty))$ satisfying $\zeta_j\equiv 1$
  in $[0,l]$ for $l\in\N$, by means of a straightforward extraction procedure we obtain $(\eps_j)_{j\in\N} \subset (0,1)$
  and $u=(u_1,...,u_N):\Omega\times (0,\infty)\to\R^N$ such that for all $\iin$ we have
  $u_0\ge 0$ and $\uie \to u_i$ a.e.~in $\Omega\times (0,\infty)$ as $\eps=\eps_j\searrow 0$.
  Since for each $\iin$ and all $T>0$ we know from Lemma \ref{lem2} that $(\uie)_{\eps\in (0,1)}$ is bounded in
  $L^{m_i+1}(\Omega\times (0,T))$, and that hence $(\uie^p)_{\eps\in (0,1)}$ is uniformly integrable over $\Omega\times (0,T)$
  for all $p\in [1,m_i+1)$, by reflexivity of $L^{m_i+1}(\Omega\times (0,T))$ and the Vitali convergence theorem
  we readily infer that on passing to a further subsequence if necessary we can also achieve simultaneous validity of
  (\ref{7.2}) amd (\ref{7.3}).
\qed
Our next goal is to show that the functions just constructed acutally form a solution. We begin by confirming that they
enjoy a renormalized supersolution property in the style of Definition \ref{defi5}. 
The most crucial ingredient in our verification of this -- and actually the reason for dealing with \textit{super}solutions -- becomes apparent in \eqref{8.4}, which is enlisted to control the integral involving the gradient (of $P_i^{(1)}(u)$) from above by means of lower semicontinuity.
\begin{lem}\label{lem8}
Let $u_1,...,u_N$ be as given by Lemma \ref{lem7}.
  Then $u=(u_1,...,u_N)$ forms a renormalized supersolution of (\ref{0}) in the sense of Definition \ref{defi5}.
\end{lem}
\proof
  We fix $\iin$ and a nonincreasing concave $\rho\in C^\infty([0,\infty))$ such that $\rho' \in C_0^\infty([0,\infty))$,
  and use (\ref{0eps}) to see that for all nonnegative $\varphi\in C_0^\infty(\bom\times [0,\infty))$ such that
  $\frac{\partial\varphi}{\partial\nu}=0$ on $\pO\times (0,\infty)$, we have
  \bea{8.1}
	& & \hspace*{-20mm}
	- \int_0^\infty \io \rho(\uie) \varphi_t
	- \io \rho(u_{0i\eps}) \varphi(\cdot,0) 
	= + \int_0^\infty \io \partial_t \rho(\uie) \varphi \nn\\
	&=& \int_0^\infty \io \rho'(\uie) \varphi \cdot \bigg\{
	d_i m_i \nabla \cdot \Big( (\uie+\eps)^{m_i-1} \nabla \uie \Big)
	+ \frac{\fiue}{1+\eps\sum_{j=1}^N |\fjue|} \bigg\} \nn\\
	&=& -d_i m_i \int_0^\infty \io (\uie+\eps)^{m_i-1} \rho''(\uie) |\nabla\uie|^2 \varphi 
	- d_i m_i \int_0^\infty \io (\uie+\eps)^{m_i-1} \rho'(\uie) \nabla\uie\cdot\nabla\varphi \nn\\
	& & + \int_0^\infty \io \rho'(\uie) \cdot \frac{\fiue}{1+\eps\sum_{j=1}^N |\fjue|} \cdot \varphi \nn\\[2mm]
	&=& - d_i m_i \int_0^\infty \io |\nabla P_{i\eps}^{(1)} (\uie)|^2 \varphi
	+ d_i m_i \int_0^\infty \io P_{i\eps}^{(2)} (\uie) \Delta\varphi \nn\\
	& & + \int_0^\infty \io \rho'(\uie) \cdot \frac{\fiue}{1+\eps\sum_{j=1}^N |\fjue|} \cdot \varphi
	\qquad \mbox{for all } \eps\in (0,1),
  \eea
  where we have set
  \be{8.2}
	P_{i\eps}^{(1)}(s):=\int_0^s (\sigma+\eps)^\frac{m_i-1}{2} \sqrt{\rho''(\sigma)} d\sigma
	\quad \mbox{and} \quad
	P_{i\eps}^{(2)}(s):=\int_0^s (\sigma+\eps)^{m_i-1} \rho'(\sigma) d\sigma,
	\qquad s\ge 0,
  \ee
  for $\iin$ and $\eps\in (0,1)$.
  Here we note that if we take $M>0$ and $T>0$ large enough fulfilling $\supp \rho' \subset [0,M]$ and 
  $\supp\varphi \subset \bom\times [0,T]$, then for all $\eps\in (0,1)$,
  \bas
	\int_0^\infty \io\Big| \sqrt{\varphi} \nabla P_{i\eps}^{(1)} (\uie) \Big|^2
	\le \|\rho''\|_{L^\infty((0,\infty))} \|\varphi\|_{L^\infty(\Omega\times (0,\infty))}
	\int_0^T \io \chi_{\{\uie\le M\}} (\uie+\eps)^{m_i-1} |\nabla\uie|^2,
  \eas
  whence again employing Lemma \ref{lem3} we infer that with $(\eps_j)_{j\in\N}$ as provided by Lemma \ref{lem7},
  we can find a subsequence $(\eps_{j_k})_{k\in\N}$ such that
  \be{8.3}
	\sqrt{\varphi} \nabla P_{i\eps}^{(1)} (\uie)
	\wto z
	\quad \mbox{in } L^2(\Omega\times (0,\infty);\R^n)
	\qquad \mbox{as } \eps=\eps_{j_k} \searrow 0
  \ee
  for some $z\in L^2(\Omega\times (0,\infty);\R^n)$.
  On the other hand, since (\ref{8.2}) entails that $P_{i\eps}^{(1)} \to P_i^{(1)}$ in $L^\infty_{loc}([0,\infty))$ as 
  $\eps\searrow 0$, with $P_i^{(1)}$ given by (\ref{5.1}), and since moreover
  \bea{8.33}
	|P_{i\eps}^{(1)}(s)|
	&\le& \|\rho''\|_{L^\infty((0,\infty))}^\frac{1}{2}
	\int_0^M (\sigma+\eps)^\frac{m_i-1}{2} d\sigma \nn\\
	&=& \|\rho''\|_{L^\infty((0,\infty))}^\frac{1}{2}
	\cdot \frac{(M+\eps)^\frac{m_i+1}{2}-\eps^\frac{m_i+1}{2}}{\frac{m_i+1}{2}} \nn\\
	&\le& \|\rho''\|_{L^\infty((0,\infty))}^\frac{1}{2}
	\cdot \frac{2(M+1)^\frac{m_i+1}{2}}{m_i+1}
	\qquad \mbox{for all $s\ge 0$ and } \eps\in (0,1),
  \eea
  from (\ref{7.1}) and the dominated convergence theorem it follows that
  \bas
	P_{i\eps}^{(1)}(\uie) \to P_i^{(1)}(u_i)
	\quad \mbox{in } L^1(\Omega\times (0,T))
	\qquad \mbox{as } \eps=\eps_j\searrow 0.
  \eas
  Therefore, a standard argument shows that in (\ref{8.3}) we must have $z=\sqrt{\varphi} \nabla P_i^{(1)}(u_i)$
  a.e.~in $\{\varphi>\delta\}$ for all $\delta>0$, and hence actually
  $\sqrt{\varphi} \nabla P_{i\eps}^{(1)} (\uie) \wto \sqrt{\varphi} \nabla P_i^{(1)}(u_i)$ 
  in $L^2_{loc}(\Omega\times (0,\infty))$ as $\eps=\eps_{j_k}\searrow 0$, so that by lower semicontinuity of the norm in
  $L^2(\Omega\times (0,\infty))$ with respect to weak convergence,
  \be{8.4}
	d_i m_i \int_0^\infty \io |\nabla P_i^{(1)}(u_i)|^2 \varphi
	\le \liminf_{\eps=\eps_{j_k}\searrow 0} \bigg\{ 
	d_i m_i \int_0^\infty \io |\nabla P_{i\eps}^{(1)}(\uie)|^2 \varphi \bigg\}.
  \ee
  Next addressing the integrals in (\ref{8.1}) exclusively containing zero-order expressions with respect to $\uie$, we first 
  observe that clearly
  \bas
	|\rho(s)| \le \|\rho'\|_{L^\infty((0,\infty))} \cdot M +\rho(0)
	\qquad \mbox{for all } s\ge 0,
  \eas
  and that furthermore, by (\ref{8.2}), similarly to (\ref{8.33}) we can estimate
  \bas
	|P_{i\eps}^{(2)}(s)|
	&\le& \|\rho'\|_{L^\infty((0,\infty))} \int_0^M (\sigma+\eps)^{m_i-1} d\sigma \\
	&\le& \|\rho'\|_{L^\infty((0,\infty))} \cdot \frac{(M+1)^{m_i}}{m_i}
	\qquad \mbox{for all $s\ge 0$ and } \eps\in (0,1).
  \eas
  Therefore, three applications of the dominated convergence theorem on the basis of (\ref{7.1}) and (\ref{ie}) show that if we 
  take $P_i^{(2)}$ from (\ref{5.1}) then 
  \be{8.5}
	\int_0^\infty \io \rho(\uie) \varphi_t
	\to \int_0^\infty \io \rho(u_i)\varphi_t
  \ee
  and 
  \be{8.6}
	\io \rho(u_{0i\eps}) \varphi(\cdot,0)
	\to \io \rho(u_{0i}) \varphi(\cdot,0)
  \ee
  as well as
  \be{8.7}
	d_i m_i \int_0^\infty \io P_{i\eps}^{(2)}(\uie) \Delta\varphi
	\to d_i m_i \int_0^\infty \io P_i^{(2)}(u_i) \Delta\varphi
  \ee
  as $\eps=\eps_j\searrow 0$, the latter because in addition obviously $P_{i\eps}^{(2)} \to P_i^{(2)}$
  in $L^\infty_{loc}([0,\infty))$ as $\eps\searrow 0$.\abs
  Finally, in the crucial rightmost summand in (\ref{8.1}) containing the reactive contribution, we once more rewrite
  $f_i=|f_i|-2(f_i)_-$ and note that fixing any $\delta>0$ such that $(1+\delta)\beta_j \le m_j+1$ for all
  $j\in\{1,...,N\} \setminus \{i\}$, again relying on (\ref{f3}) we can estimate
  \bas
 	& & \hspace*{-20mm}
 	\int_0^T \io \bigg| \rho'(\uie) \cdot \frac{(\fiue)_{-}}{1+\eps\sum_{j=1}^N |\fjue|} \cdot\varphi \bigg|^{1+\delta}
 	\nn\\
 	&\le& \int_0^T \io |\rho'(\uie)|^{1+\delta} \cdot 
 	\bigg\{ \phi_i(\uie) \bigg(\sum_{j\ne i} u_{j\eps}^{\beta_j} +1\bigg) \bigg\}^{1+\delta} \cdot \varphi^{1+\delta} \\
 	&\le& \Big\{ \|\rho'\|_{L^\infty((0,\infty))} \cdot \|\phi_i\|_{L^\infty((0,M))} \cdot
 		\|\varphi\|_{L^\infty(\Omega\times (0,\infty))} \Big\}^{1+\delta}
 	\cdot N^{1+\delta} \cdot\bigg(
 	\sum_{j\ne i} \int_0^T \io u_{j\eps}^{(1+\delta)\beta_j} +|\Omega|T\bigg)
  \eas
  for all $\eps\in (0,1)$.
  In view of Lemma \ref{lem2}, by positivity of $\delta$ this implies uniform integrability of
  $\Big(\rho'(\uie) \cdot \frac{(\fiue)_-}{1+\eps\sum_{j=1}^N |\fjue|} \cdot\varphi\Big)_{\eps\in (0,1)}$ 
  over $\Omega\times (0,T)$ and hence entails, when combined with (\ref{7.1}), that
  \be{8.8}
	- 2 \int_0^\infty \io \rho'(\uie) \cdot \frac{(\fiue)_-}{1+\eps\sum_{j=1}^N |\fjue|} \cdot \varphi
	\to -2 \int_0^\infty \io \rho'(u_i) (f_i(u_1,...,u_N))_- \varphi
  \ee
  as $\eps=\eps_j\searrow 0$.
  Since apart from that, by nonnegativity of both $-\rho'$ and $\varphi$ we can invoke Fatou's lemma to see that, again 
  thanks to (\ref{7.1}),
  \be{8.9}
	\hspace*{-5mm}
	- \int_0^\infty \io \rho'(u_i) |f_i(u_1,...,u_N)| \varphi
	\le \liminf_{\eps=\eps_j\searrow 0} \bigg\{ 
	- \int_0^\infty \io \rho'(\uie) \cdot \frac{|\fiue|}{1+\eps\sum_{j=1}^N |\fjue|} \cdot \varphi \bigg\},
  \ee
  upon collecting (\ref{8.4})-(\ref{8.9}) we altogether conclude from (\ref{8.1}) that
  \bas
	& & \hspace*{-14mm}
	d_i m_i \int_0^\infty \io |\nabla P_i^{(1)} (u_i) |^2 \varphi
	- \int_0^\infty \io \rho'(u_i) |f_i(u_1,...,u_N)| \varphi \\
	&\le& \liminf_{\eps=\eps_{j_k}\searrow 0}
	\bigg\{ d_i m_i \int_0^\infty \io \Big|\nabla P_{i\eps}^{(1)}(\uie)|^2 \varphi
	- \int_0^\infty \io \rho'(\uie) \cdot \frac{|\fiue|}{1+\eps\sum_{j=1}^N |\fjue|} \cdot \varphi \bigg\} \\[2mm]
	&=& \liminf_{\eps=\eps_{j_k}\searrow 0}
	\bigg\{ \int_0^\infty \io \rho(\uie) \varphi_t
	+ \io \rho(u_{0i\eps})\varphi(\cdot,0) \\
	& & \hspace*{16mm}
	+ d_i m_i \int_0^\infty \io P_{i\eps}^{(2)} (\uie) \Delta \varphi
	- 2 \int_0^\infty \io \rho'(\uie) \cdot \frac{(\fiue)_-}{1+\eps\sum_{j=1}^N |\fjue|} \cdot \varphi \bigg\} \\[2mm]
	&=& 
	\int_0^\infty \io \rho(u_i) \varphi_t
	+ \io \rho(u_{0i})\varphi(\cdot,0) \\
	& & \hspace*{16mm}
	+ d_i m_i \int_0^\infty \io P_i^{(2)} (u_i) \Delta \varphi
	- 2 \int_0^\infty \io \rho'(u_i) (f_i(u_1,...,.u_N))_- \varphi,
  \eas
  which is equivalent to the desired inequality (\ref{5.4}).
  The integrablity requirements in (\ref{5.23}) are evident by-products of the above considerations.
\qed
But also the subsolution property encoded in \eqref{6.2} is fulfilled:
\begin{lem}\label{lem89}
  The function $u=(u_1,...,u_N)$ from Lemma \ref{lem7} 
  satisfies \eqref{6.1} and \eqref{6.2} of Definition~\ref{defi6}.
%
\end{lem}
\proof
  According to (\ref{7.2}), we can pick a null set $N\subset (0,\infty)$ such that with $(\eps_j)_{j\in\N}$ taken from
  Lemma \ref{lem7}, for each $t\in (0,\infty)\setminus N$ we have $\uie(\cdot,t)\to u_i(\cdot,t)$ in $L^1(\Omega)$ for 
  all $\iin$ and hence
  \be{89.1}
	\io \sum_{i=1}^N a_i \uie(\cdot,t) \to \io \sum_{i=1}^N a_i u_i(\cdot,t)
  \ee
  as $\eps=\eps_j\searrow 0$, where $a\in(0,\infty)^N$ is taken from \eqref{f1}.
  We next let $F_\eps(s_1,...,s_N):=\sum_{i=1}^N \frac{a_i  f_i(s_1,...,s_N)}{1+\eps\sum_{j=1}^N |f_j(s_1,...,s_N)|}$
  and $F(s_1,...,s_N):=\sum_{i=1}^N a_i  f_i(s_1,...,s_N)$ for $(s_1,...,s_N)\in \R^N$ and $\eps\in (0,1)$,
  and then obtain from (\ref{f1}) that
  \bas
	\Big(F_\eps(u_{1\eps},...,u_{N\eps})\Big)_+
	\le K \sum_{i=1}^N a_i \uie + K
	\quad \mbox{in } \Omega\times (0,\infty)
	\qquad \mbox{for all } \eps\in (0,1).
  \eas
  In view of (\ref{7.2}) applied to $p:=1$, a version of the dominated convergence theorem thus ensures that
  \be{89.2}
	\int_0^t \io \Big(F_\eps(u_{1\eps},...,u_{N\eps})\Big)_+
	\to \int_0^t \io F_+(u_1,...,u_N)
	\quad \mbox{for all } t>0
	\qquad \mbox{as } \eps=\eps_j\searrow 0,
  \ee
  because clearly $\Big(F_\eps(u_{1\eps},...,u_{N\eps})\Big)_+ \to F_+(u_1,...,u_N)$ a.e.~in $\Omega\times (0,\infty)$
  as $\eps=\eps_j\searrow 0$ by (\ref{7.1}).
  Since from (\ref{0eps}) we know that for all $t>0$ and $\eps\in (0,1)$ we have
  \bas
	\io \sum_{i=1}^N a_i \uie(\cdot,t)
	+ \int_0^t \io 	\Big(F_\eps(u_{1\eps},...,u_{N\eps})\Big)_-
	= \io \sum_{i=1}^N a_i u_{0i\eps}
	+ \int_0^t \io 	\Big(F_\eps(u_{1\eps},...,u_{N\eps})\Big)_+,
  \eas
  where by (\ref{7.1}) also $\Big(F_\eps(u_{1\eps},...,u_{N\eps})\Big)_- \to F_-(u_1,...,u_N)$ a.e.~in $\Omega\times (0,\infty)$
  as $\eps=\eps_j\searrow 0$, and where $\io \sum_{i=1}^N a_i u_{0i\eps} \to \io \sum_{i=1}^N a_i u_{0i}$ as $\eps\searrow 0$ due to
  (\ref{ie}), invoking Fatou's lemma we infer by means of (\ref{89.1}) and (\ref{89.2}) that
  \bas
	\io \sum_{i=1}^N a_i  u_i(\cdot,t)
	+ \int_0^t \io F_-(u_1,...,u_N)
	\le \io \sum_{i=1}^N a_i  u_{0i}
	+ \int_0^t \io F_+(u_1,...,u_N)
	\qquad \mbox{for all } t\in (0,\infty)\setminus N.
  \eas
  For any such $t$, this firstly implies that $F(u_1,\ldots,u_N)$ belongs to $L^1(\Omega\times (0,t))$, and secondly entails that
  (\ref{6.2}) holds.
\qed
The previous two lemmata already demonstrate that $u$ is a generalized solution in the sense of Definition \ref{defi6}:\abs
\proofc of Theorem \ref{theo9}. \quad
  We take $u_1,...,u_N$ as given by Lemma \ref{lem7} and then only need to combine Lemma \ref{lem8} with Lemma \ref{lem89}.
\qed

{
\footnotesize
\setlength{\parskip}{0pt}
\setlength{\itemsep}{0pt plus 0.3ex}
\def\cprime{$'$}

}

\end{document}